\documentclass[12pt]{amsart}

\usepackage[english]{babel}
\usepackage[letterpaper,top=3cm,bottom=3cm,left=3.25cm,right=3.25cm,marginparwidth=1.75cm]{geometry}
\usepackage{amsmath}
\numberwithin{equation}{section}
\usepackage{amsthm}
\usepackage{amssymb}
\usepackage{mathtools}
\usepackage{bm}
\usepackage{dsfont}
\usepackage{array}
\usepackage{graphicx}
\usepackage{float}
\usepackage{flafter}
\usepackage{booktabs}
\usepackage{siunitx}
\usepackage{placeins}
\usepackage[dvipsnames]{xcolor}
\usepackage{tikz}
\usetikzlibrary{arrows.meta,positioning,calc,backgrounds,fit}
\usepackage[colorlinks=true,linkcolor=blue,citecolor=green,urlcolor=green]{hyperref}

\newcommand{\PaperTitle}{Selective boundary condition reduction via learned error gating}

\newcommand{\PaperAbstract}{%
  Parametric PDEs can admit different boundary conditions with different
  accuracy and computational cost. We introduce a framework for learning when
  one reduced boundary condition can replace another: paired solutions train a
  neural network to estimate the resulting domain and boundary errors, and the
  simpler condition is used only when both predicted errors meet prescribed
  tolerances.
  We focus on singular limits in applications, in which a
  stiff Robin or nonlinear boundary law is replaced by its limiting Dirichlet
  form. We evaluate the method on a galvanic corrosion problem and other
  nonlinear stationary and evolution problems.%
}

\newcommand{\PaperKeywords}{Parametric PDEs, singular boundary limits, Robin
boundary conditions, error estimation, model selection}

\definecolor{legBothReject}{HTML}{F2F2F2}
\definecolor{legBothAccept}{HTML}{1B9E77}
\definecolor{legUnsafe}{HTML}{FF0000}
\definecolor{legMissed}{HTML}{0000FF}

\newtheorem{theoremp}{Theorem}[section]
\newtheorem{propositionp}[theoremp]{Proposition}

\theoremstyle{definition}
\newtheorem{assumptionp}[theoremp]{Assumption}
\newenvironment{assumption}{\vspace{5pt}\begin{assumptionp}}{\end{assumptionp}\vspace{5pt}}

\makeatletter
\renewcommand{\@defaultbiblabelstyle}[1]{[#1]}
\newtheoremstyle{boldremark}
  {5pt}
  {5pt}
  {\normalfont}
  {}
  {\bfseries}
  {.}
  { }
  {}
\makeatother
\theoremstyle{boldremark}

\newcolumntype{L}[1]{>{\raggedright\arraybackslash}p{#1}}
\newcommand{\tablecaptiongap}{\vspace{6pt}}

\title[Selective Boundary Condition Reduction]{\MakeUppercase{\PaperTitle}}
\author[D. Fern\'andez, D. Penk and D. Riedelbauch]{%
Daniel Fern\'andez\textsuperscript{$\ast$},
Dominik Penk\textsuperscript{$\dagger$} and
Dominik Riedelbauch\textsuperscript{$\dagger$}}
\thanks{%
\textsuperscript{$\ast$}Chair for Dynamics, Control, Machine Learning, and
Numerics (Alexander von Humboldt Professorship), Department of Mathematics,
Friedrich-Alexander-Universit\"at Erlangen-N\"urnberg, Cauerstra{\ss}e,
91058 Erlangen, Germany. \emph{Email:}
\texttt{daniel.fernandez@fau.de}.\newline
\textsuperscript{$\dagger$}Schaeffler Technologies AG \& Co. KG,
Industriestra{\ss}e 1--3, Herzogenaurach, Germany. \emph{Emails:}
\texttt{penkdmi@schaeffler.com}, \texttt{riededmi@schaeffler.com}.}
\date{}

\begin{document}

\begin{abstract}
\PaperAbstract
\end{abstract}

\maketitle

\providecommand{\paperbeginstationaryfigure}{\begin{figure}[htbp]}
\providecommand{\paperbeginstationarytable}{\begin{table}[htbp]}
\providecommand{\paperbeginmethodfigure}{\begin{figure*}[!t]}
\providecommand{\paperendmethodfigure}{\end{figure*}}
\providecommand{\paperbegincorrosiontable}{\begin{table}[htbp]}
\providecommand{\paperendcorrosiontable}{\end{table}}
\providecommand{\paperbegincorrosionfigure}{\begin{figure}[htbp]}
\providecommand{\paperendcorrosionfigure}{\end{figure}}
\providecommand{\papermethodbarrier}{\FloatBarrier}
\providecommand{\papercorrosionbarrier}{\FloatBarrier}

\section{Introduction}\label{sec:introduction}

Boundary conditions govern the exchange of heat, mass, or charge with the
surroundings. Different conditions can describe the same exchange process at
different fidelity and cost, motivating us to learn whether one can replace
another for each parameter instance. We focus on singular limits in
applications. Under strong transfer or reaction, a stiff Robin or
nonlinear law drives the boundary value toward a target; electrochemical
interfaces are one example \cite{Bazant2013}. The full law can be expensive to
solve even though its singular limit reduces to a Dirichlet condition.

Classical analyses establish convergence to the Dirichlet problem
\cite{Costabel1996,Auchmuty2018,Amrouche2020}, but not whether the limiting law
is accurate enough at a given finite parameter value. That decision also
depends on the forcing, coefficients, geometry, and error norm. In particular,
a small domain error can coexist with a much larger boundary error that matters
for coupling.

In contrast to PINNs, neural operators, and other neural surrogates, which
usually fail in singular regimes
\cite{Gie2024,Takamoto2022PDEBench,Lanthaler2023Nonlinear,
Li2024ComponentFNO}, our selector avoids reconstructing the full field. It
learns two scalar discrepancies, while a PDE solver returns the field and
enforces the selected law. The learned component chooses which PDE model is
solved.
The full law remains the fallback, and new tolerances require no
retraining. By targeting the decision quantities while retaining solver
accuracy, the method can outperform direct field emulation in data efficiency
and reliability when singular or localized features control the tolerance
decision.

\subsection{Contributions}

The paper makes three main contributions:
\begin{itemize}
  \item \textbf{Selection supported by PDE solvers in singular regimes.}
  Predicting two errors avoids full field reconstruction in singular regimes.
  The solution is computed by a PDE solver that enforces the selected boundary
  law. The full law remains the fallback. Tolerances can change without
  retraining, and split conformal calibration can make both error estimates
  conservative.

  \item \textbf{Evaluation on stationary and evolutionary PDEs.}
  Three benchmarks cover harmonic galvanic corrosion, stationary nonlinear
  transfer, and nonlinear evolution. In the tested stationary and evolutionary
  examples, the predictor identifies when the simpler boundary law meets the
  prescribed tolerances, with few observed selection errors.

  \item \textbf{Speedup of the complete online policy.}
  The stationary benchmark achieves about $133\times$ speedup when the limit
  model is selected and $8.3\times$ for the complete policy. Timings include
  feature construction, inference, selected solves, and full fallbacks; offline
  preparation is excluded.
\end{itemize}

\subsection{Related work}

Analyses of limits from Robin to Dirichlet conditions give convergence rates
and uniform estimates \cite{Costabel1996,Auchmuty2018,Amrouche2020}, including
results for contact impedance and nonlinear jumps in boundary data
\cite{Darde2016,Fernandez2026SingularLimit};
robust discretizations across boundary regimes appear in \cite{Juntunen2009}.
This literature establishes or discretizes the limit; we learn whether it meets
separate domain and trace tolerances at a finite parameter instance.

Neural surrogate models, including PINNs and neural operators
\cite{Kovachki2023}, are often unreliable near singular phenomena such as
shocks, boundary layers, and singular limits
\cite{Gie2024,Takamoto2022PDEBench,LiuSchiaffini2024Localized,
Liu2024SpectralBias,Lanthaler2023Nonlinear,Li2024ComponentFNO}. Our selector
avoids full field reconstruction and learns only the errors needed to choose the
boundary condition.

Adaptivity based on quantities of interest, certified hierarchies, learned
error estimates, and hybrid switching also address model choice
\cite{Oden2001,Haasdonk2023,Freno2019,Fritzen2019,Riffaud2025}. Our
specialization holds the state space, interior operator, and discretization
fixed, changing only the boundary law. The rule relates to selective prediction
\cite{Chow1970,Geifman2017}; split conformal calibration gives marginal coverage
under exchangeability \cite{Lei2018}.

\section{Mathematical setting and structure of the singular limit}\label{sec:model}

\subsection{A preliminary elliptic example}

Let $\Omega\subset\mathbb R^d$ be a bounded Lipschitz domain, let
$\Gamma\subset\partial\Omega$ be relatively open, and let $T$ denote the trace
on $\Gamma$. For $f\in L^2(\Omega)$, target $u_D\in L^2(\Gamma)$, and
$\kappa>0$, the Robin problem
\begin{equation}
  u_\kappa
    =\operatorname*{argmin}_{v\in H^1(\Omega)}
      \left\{\frac12\int_\Omega(|\nabla v|^2+|v|^2)\,dx
      -\int_\Omega fv\,dx
      +\frac{1}{2\kappa}\int_\Gamma |Tv-u_D|^2\,ds\right\}
  \label{eq:preliminary-penalty-energy}
\end{equation}
has Euler equation $-\Delta u_\kappa+u_\kappa=f$ and boundary law
$\partial_nu_\kappa+\kappa^{-1}(u_\kappa-u_D)=0$ on $\Gamma$ (with homogeneous
Neumann data elsewhere). As $\kappa\downarrow0$, the last term in
\eqref{eq:preliminary-penalty-energy} enforces $Tu=u_D$. If
$\{v\in H^1(\Omega):Tv=u_D\}$ is nonempty, the minimizers converge to the
corresponding problem constrained by Dirichlet data. For a nonzero jump in
$u_D$, that set may be empty and the limit must instead be interpreted in a
weaker trace class; this is the situation in Experiment~1. Standard analyses of Robin limits
make these statements precise \cite{Costabel1996,Auchmuty2018}.

\subsection{Full and limit models}

We compare a full solution $u_{\mathrm{full}}$ with a limit solution
$u_{\mathrm{lim}}$ on a domain $\Omega$. Away from a designated boundary part
$\Gamma$, both satisfy
\begin{equation}
  \mathcal L u=f \quad\text{in }\Omega,
  \label{eq:common-interior-problem}
\end{equation}
and the same conditions on $\partial\Omega\setminus\Gamma$. The operator,
source, boundary response, geometry, and prescribed data may vary from case to
case. We suppress that dependence in the notation; only the law on $\Gamma$
changes.

The full model uses the possibly nonlinear law
\begin{equation}
  \partial_n u_{\mathrm{full}}+\frac{1}{\kappa}
  \beta(x,u_{\mathrm{full}})=0
  \qquad\text{on }\Gamma,
  \label{eq:nonlinear_neumann_penalty}
\end{equation}
where $n$ is the outward unit normal and boundary values denote traces. The
response $\beta$ drives the trace toward a preferred value $u_D(x)$, which may
differ between boundary pieces, and smaller $\kappa>0$ means a stiffer response.
The limit model replaces
\eqref{eq:nonlinear_neumann_penalty} by
\begin{equation}
  u_{\mathrm{lim}}=u_D(x) \qquad\text{on }\Gamma.
  \label{eq:dirichlet-limit-law}
\end{equation}

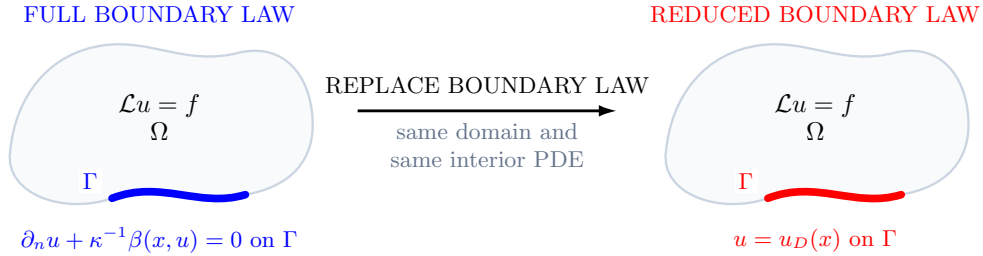
\begin{figure}[htbp]
\centering
\resizebox{0.9\linewidth}{!}{%
  \begingroup
\definecolor{bcMuted}{HTML}{64748B}
\definecolor{bcPanel}{HTML}{F8FAFC}
\definecolor{bcLine}{HTML}{CBD5E1}
\begin{tikzpicture}[
  font=\normalfont\footnotesize,
  domainshape/.style={
    draw=bcLine,
    fill=bcPanel,
    line width=0.9pt
  },
  transition/.style={
    -{Latex[length=2.5mm,width=1.6mm]},
    draw=black,
    line width=1.1pt
  },
  title/.style={font=\normalfont\scriptsize},
  heading/.style={font=\normalfont\bfseries\scriptsize},
  note/.style={font=\normalfont\scriptsize,text=bcMuted}
]

\path[domainshape]
  (-2.05,-0.35)
  .. controls (-1.98,0.42) and (-1.52,1.02) .. (-0.76,1.12)
  .. controls (-0.10,1.21) and (0.36,0.83) .. (1.10,1.02)
  .. controls (1.72,1.16) and (2.07,0.69) .. (2.10,0.12)
  .. controls (2.14,-0.45) and (1.81,-0.82) .. (1.18,-0.93)
  .. controls (0.55,-1.12) and (0.02,-0.72) .. (-0.65,-0.98)
  .. controls (-1.32,-1.21) and (-2.12,-0.92) .. cycle;
\begin{scope}[shift={(9,0)}]
  \path[domainshape]
    (-2.05,-0.35)
    .. controls (-1.98,0.42) and (-1.52,1.02) .. (-0.76,1.12)
    .. controls (-0.10,1.21) and (0.36,0.83) .. (1.10,1.02)
    .. controls (1.72,1.16) and (2.07,0.69) .. (2.10,0.12)
    .. controls (2.14,-0.45) and (1.81,-0.82) .. (1.18,-0.93)
    .. controls (0.55,-1.12) and (0.02,-0.72) .. (-0.65,-0.98)
    .. controls (-1.32,-1.21) and (-2.12,-0.92) .. cycle;
\end{scope}

\node[title,text=blue] at (0,1.55) {FULL BOUNDARY LAW};
\node[title,text=red] at (9,1.55) {REDUCED BOUNDARY LAW};

\node[text=black,align=center] at (0,0.13)
  {$\mathcal L u=f$\\[-1pt]$\Omega$};
\node[text=black,align=center] at (9,0.13)
  {$\mathcal L u=f$\\[-1pt]$\Omega$};

\draw[blue,line width=3pt,line cap=round]
  (-0.65,-0.98) .. controls (0.02,-0.72) and (0.55,-1.12) .. (1.18,-0.93);
\draw[red,line width=3pt,line cap=round]
  (8.35,-0.98) .. controls (9.02,-0.72) and (9.55,-1.12) .. (10.18,-0.93);
\node[heading,text=blue,fill=white,inner sep=1.5pt]
  at (-0.95,-0.74) {$\Gamma$};
\node[heading,text=red,fill=white,inner sep=1.5pt]
  at (8.05,-0.74) {$\Gamma$};

\node[note,text=blue,align=center] at (0,-1.55)
  {$\partial_nu+\kappa^{-1}\beta(x,u)=0$ on $\Gamma$};
\node[note,text=red,align=center] at (9,-1.55)
  {$u=u_D(x)$ on $\Gamma$};

\draw[transition] (2.72,0.22) -- (6.28,0.22);
\node[title,text=black,fill=white,inner xsep=3pt]
  at (4.5,0.58) {REPLACE BOUNDARY LAW};
% Keep the note below the shaft without an opaque background over the arrow.
\node[note,align=center,anchor=north,inner xsep=2pt,inner ysep=0pt]
  at (4.5,0.06) {same domain and\\same interior PDE};

\end{tikzpicture}
\endgroup%
}
\caption{Only the boundary law changes; the domain and interior PDE remain
fixed.}
\label{fig:boundary-condition-transition}
\end{figure}

\begin{assumption}[Standing limit assumptions]\label{ass:limit}
For the parameter range under consideration, the full and limit models are well
posed in the solution classes used below. Away from finitely many junction
points, the boundary response $\beta(x,\cdot)$ has $u_D(x)$ as its only zero in
a fixed neighborhood and is uniformly strongly monotone there. The full
solution converges to the stated limit in the domain and boundary norms used for
selection.
\end{assumption}

Under Assumption~\ref{ass:limit}, the boundary identity gives
$\beta(x,u_{\mathrm{full}})
=-\kappa\partial_n u_{\mathrm{full}}$. For a smooth target and bounded normal
derivatives, strong monotonicity directly gives
$u_{\mathrm{full}}|_\Gamma\to u_D=u_{\mathrm{lim}}|_\Gamma$, and stability
carries the convergence into the domain. At a jump in $u_D$, the normal
derivative need not remain bounded and convergence instead holds in weaker
trace norms and locally in the interior. These hypotheses are standard and, in
fact, are satisfied by all three benchmarks in
Section~\ref{sec:experiments}, with the evolutionary problem interpreted in its
natural time-dependent weak formulation. Precise results show that the rate
depends on geometry, smoothness, and the chosen norm
\cite{Costabel1996,Auchmuty2018,Fernandez2026SingularLimit}.

\subsection{Error measures and the reference rule}

Let $\|\cdot\|_\Omega$ and $\|\cdot\|_\Gamma$ denote the domain and
boundary norms, respectively. The benchmarks in Section~\ref{sec:experiments}
specify how these norms are discretized. The two relative errors are
\begin{equation}
  \begin{aligned}
  E_\Omega
    &=\frac{\|u_{\mathrm{full}}-u_{\mathrm{lim}}\|_\Omega}
            {\|u_{\mathrm{lim}}\|_\Omega},&
  E_\Gamma
    &=\frac{\|u_{\mathrm{full}}-u_{\mathrm{lim}}\|_\Gamma}
            {\|u_{\mathrm{lim}}\|_\Gamma}.
  \end{aligned}
  \label{eq:ideal-boundary-reduction-reference}
\end{equation}
$E_\Omega$ measures the global change, whereas $E_\Gamma$ exposes a discrepancy
that the domain norm can hide. If both solutions were known, the ideal rule
would use the limit model exactly when
\begin{equation}
  E_\Omega\le\varepsilon_\Omega
  \qquad\text{and}\qquad
  E_\Gamma\le\varepsilon_\Gamma.
  \label{eq:ideal-selector}
\end{equation}
A choice is \emph{unsafe} if the limit model is used while either inequality
fails. Because evaluating this rule requires the full solution, the practical
rule estimates both errors before the full solve.

\paragraph{Errors over an evolutionary trajectory.}

When one model is chosen for a complete trajectory on $0\le t\le T$, the two
errors are the largest discrepancies between the full and limit solutions over
the stored time levels:
\begin{equation}
  E_j^{\max}
  =\max_{0\le t\le T}
   \frac{\|u_{\mathrm{full}}(t)-u_{\mathrm{lim}}(t)\|_j}
        {\|u_{\mathrm{full}}(t)\|_j},
  \qquad j\in\{\Omega,\Gamma\}.
  \label{eq:evolution-trajectory-error}
\end{equation}
The evolutionary benchmark estimates both quantities before advancing the
selected trajectory and accepts the limit model only when both predicted
maxima satisfy their tolerances.

The stationary studies normalize with the limit solution, consistent with
their estimates for the singular limit. The evolutionary errors use the full
trajectory as their reference and protect every denominator by a positive
numerical floor, which is inactive in the reported cases.

\section{Method for selecting boundary laws}\label{sec:method}

The selection procedure has an offline and an online stage. Offline, paired
full and limit solves are used to learn the domain and boundary errors; online,
the predicted errors are compared with tolerances chosen by the user to select
the boundary law. We next explain the choice of learning target and how the
paired data are organized.

\subsection*{Choosing the error estimate and applying the rule}

Paired solutions can support several learning targets: a classifier tied to one
tolerance pair, the two scalar errors, a correction field, or the full
solution. We predict the two errors. This keeps the output dimension
independent of the mesh, allows tolerances to change without retraining, and
ensures that the returned field still comes from a PDE solver. By contrast,
field surrogates have many more outputs and are more costly to predict and
verify.

Within this formulation, we use neural networks primarily because the errors depend nonlinearly on several interacting PDE parameters.

\paperbeginmethodfigure
\centering
\begingroup%
\definecolor{fdInk}{HTML}{1A1A1A}%
\definecolor{fdMute}{HTML}{6E7681}%
\definecolor{fdRule}{HTML}{AEB5BF}%
\definecolor{fdBand}{HTML}{EBECEF}%
\colorlet{fdBlue}{blue}%
\colorlet{fdRed}{red}%
\begin{tikzpicture}[
  x=1cm,y=1cm,
  font=\normalfont\footnotesize,
  box/.style={draw=fdRule,line width=0.5pt,fill=white,align=center,
              text=fdInk,inner xsep=4pt,inner ysep=4pt},
  fullsolve/.style={box,fill=fdBlue!4!white,draw=fdBlue!35},
  limsolve/.style={box,fill=fdRed!3!white,draw=fdRed!35},
  solvebar/.style={line width=2.1pt,line cap=butt},
  est/.style={box,draw=fdInk,line width=0.7pt},
  dec/.style={box,rounded corners=7pt,draw=fdInk,line width=0.9pt},
  flow/.style={-{Latex[length=2.1mm,width=1.4mm]},draw=fdInk,
               line width=0.75pt,rounded corners=1.5pt},
  hair/.style={-{Latex[length=1.9mm,width=1.25mm]},draw=fdMute,
               line width=0.5pt,rounded corners=1.5pt},
  cond/.style={-{Latex[length=1.9mm,width=1.25mm]},draw=fdRed,
               line width=0.6pt,rounded corners=1.5pt,
               dash pattern=on 2.4pt off 1.5pt},
  note/.style={font=\normalfont\scriptsize,text=fdMute,align=center,inner sep=1pt},
  ptitle/.style={font=\normalfont\footnotesize,text=fdInk,anchor=south west,inner sep=0pt},
  pscope/.style={font=\normalfont\scriptsize,text=fdMute,anchor=south east,inner sep=0pt}
]

\node[box,text width=3.02cm]       (a1)  at (1.680,0)     {parameter samples\\[1pt]$i=1,\dots,N$};
\node[fullsolve,text width=2.65cm] (a2f) at (5.050,0.44)   {full solve $u_{\mathrm{full},i}$};
\node[limsolve,text width=2.65cm]  (a2l) at (5.050,-0.44)  {limit solve $u_{\mathrm{lim},i}$};
\node[box,text width=2.15cm]       (a3)  at (8.150,0)      {paired errors\\[1pt]$E_{\Omega,i},\;E_{\Gamma,i}$};
\node[est,text width=3.70cm]       (a4)  at (12.100,0)
   {fit estimator $p_\theta$\\[1pt]$\varphi_i\mapsto(\widehat m_\Omega,\widehat m_\Gamma)$\\[1pt]
    {\scriptsize\textcolor{fdMute}{$c\ge1$ from separate cases}}};

\draw[solvebar,fdBlue] ([xshift=1.15pt]a2f.north west) -- ([xshift=1.15pt]a2f.south west);
\draw[solvebar,fdRed]  ([xshift=1.15pt]a2l.north west) -- ([xshift=1.15pt]a2l.south west);

\node[note] (na2) at (5.050,-0.98)  {two solves per sample};
\node[note] (na3) at (8.150,-0.98)  {$m_j=\log_{10}E_j$};

\draw[flow] (a1.east)  -- ++(0.115,0) |- (a2f.west);
\draw[flow] (a1.east)  -- ++(0.115,0) |- (a2l.west);
\draw[flow] (a2f.east) -- ++(0.20,0) |- (a3.west);
\draw[flow] (a2l.east) -- ++(0.20,0) |- (a3.west);
\draw[flow] (a3.east)  -- (a4.west);

\draw[hair] (a1.north) -- (1.680,1.02) -- (12.100,1.02) -- (a4.north);
\node[note,fill=white,inner sep=1.2pt] (nphi) at (7.100,1.02) {features $\varphi_i$};

\node[ptitle] (ta) at (0.005,1.40)  {\textbf{(a)}\; Offline preparation};
\node[pscope] (sa) at (14.091,1.40) {$N$ paired training cases};

\begin{scope}[shift={(0,-3.78)}]
\node[box,text width=3.02cm] (b1) at (1.680,0) {new case\\[1pt]parameters $\mu$};
\node[est,text width=2.65cm] (b2) at (5.050,0)
   {evaluate $p_\theta(\varphi)$\\[1pt]$\widehat E_j=10^{\widehat m_j}$};
\node[dec,text width=2.15cm] (b3) at (8.150,0)
   {$\left\{\begin{array}{@{}l@{}}c\widehat E_\Omega\le\varepsilon_\Omega\\[1pt]
                                  c\widehat E_\Gamma\le\varepsilon_\Gamma\end{array}\right.$};
\node[note] (nb3a) at (8.150,1.10) {tolerance test};

\node[box,text width=3.45cm] (bF) at (12.100,0.80)
   {{\scriptsize\textcolor{fdMute}{either bound fails}}\\[1pt]\textcolor{fdBlue}{full law}\\[1pt]
    {\scriptsize added full solve $\to u_{\mathrm{full}}$}};
\node[box,text width=3.45cm] (bL) at (12.100,-0.80)
   {{\scriptsize\textcolor{fdMute}{both bounds hold}}\\[1pt]\textcolor{fdRed}{limit law}\\[1pt]
    {\scriptsize no full solve $\to u_{\mathrm{lim}}$}};

\node[limsolve,text width=3.02cm] (b0) at (1.680,-1.74)
   {limit solve $u_{\mathrm{lim}}$\\[1pt]
    {\scriptsize\textcolor{fdMute}{only if $\varphi$ uses $u_{\mathrm{lim}}$}}};

\draw[solvebar,fdBlue] ([xshift=1.15pt]bF.north west) -- ([xshift=1.15pt]bF.south west);
\draw[solvebar,fdRed] ([xshift=1.15pt]bL.north west) -- ([xshift=1.15pt]bL.south west);
\draw[solvebar,fdRed]  ([xshift=1.15pt]b0.north west) -- ([xshift=1.15pt]b0.south west);

\draw[flow] (b1.east) -- (b2.west);
\draw[flow] (b2.east) -- (b3.west);
\draw[flow] (b3.east) -- (9.850,0) -- (9.850,0.80)  -- (bF.west);
\draw[flow] (b3.east) -- (9.850,0) -- (9.850,-0.80) -- (bL.west);

\draw[cond] (b1.south) -- (b0.north);
\draw[cond,-] (b0.east) -- (3.470,-1.74) -- (3.470,0);
\fill[fdInk] (3.470,0) circle (0.043);
\node[font=\normalfont\scriptsize,text=fdRed,anchor=west,inner sep=1pt]
   at (3.560,-1.50) {$u_{\mathrm{lim}}$ summaries};

\draw[cond] (b0.south) -- (1.680,-2.34) -- (12.100,-2.34) -- (bL.south);
\node[font=\normalfont\scriptsize,text=fdRed,fill=white,inner sep=1.2pt]
   (nreuse) at (7.100,-2.34) {$u_{\mathrm{lim}}$ reused};

\node[ptitle] (tb) at (0.005,1.70)  {\textbf{(b)}\; Online selection};
\node[pscope] (sb) at (14.091,1.70) {one decision per case};
\end{scope}

\begin{scope}[on background layer]
  \node[fill=fdBand,inner xsep=6pt,inner ysep=6pt,
        fit=(a1)(a2f)(a2l)(a3)(a4)(na2)(na3)(nphi)(ta)(sa)] (pa) {};
  \node[inner xsep=6pt,inner ysep=6pt,
        fit=(b1)(b2)(b3)(bF)(bL)(b0)(nb3a)(nreuse)(tb)(sb)] (pb) {};
\end{scope}

\draw[flow,dash pattern=on 2.6pt off 1.8pt]
   (a4.south) -- (12.100,-1.62) -- (5.050,-1.62) -- (b2.north);
\node[font=\normalfont\scriptsize,text=fdInk,fill=white,inner sep=1.3pt]
   at (7.800,-1.62) {fitted $p_\theta$ and $c$};

\path[use as bounding box] (pa.north west) rectangle (pb.south east);

\end{tikzpicture}%
\endgroup%
\caption{Initial preparation and repeated selection between boundary laws.
\textbf{Offline}, paired full and limit solutions are used to train the error estimator.
\textbf{Online}, the selector returns the solution from the chosen model.}
\label{fig:method}
\paperendmethodfigure

Inputs containing only parameters require only the selected solve. Inputs
derived from the limit solution require that solve before selection, so savings
occur only when avoided full solves outweigh this overhead.

\subsection*{Paired solutions and errors to be estimated}

For each stationary training case, we solve both models and record the errors
in \eqref{eq:ideal-boundary-reduction-reference}. For evolutionary cases, we
record the trajectory errors in \eqref{eq:evolution-trajectory-error}.

Let $j\in\{\Omega,\Gamma\}$ denote either error type, corresponding respectively
to the interior domain and the boundary. Let $m_j:=\log_{10}E_j$, and let
$p_\theta$ denote the fitted error estimator with input vector $\varphi$. It
returns
\begin{equation}
  (\widehat m_\Omega,\widehat m_\Gamma)
  =p_\theta(\varphi),
  \qquad \widehat E_j:=10^{\widehat m_j}.
  \label{eq:log-error-estimator}
\end{equation}
Inputs contain parameters and, when used, inexpensive summaries of the limit
solution.
The practical rule chooses the limit model exactly when
\begin{equation}
  \widehat E_\Omega\le\varepsilon_\Omega
  \qquad\text{and}\qquad
  \widehat E_\Gamma\le\varepsilon_\Gamma.
  \label{eq:threshold-decision}
\end{equation}

\paragraph{Optional conservative calibration.}
To reduce unsafe selections, we multiply both predicted errors by a common
factor $c\ge1$ learned from separate calibration cases. For example, $c=1.2$
raises both estimates by 20\% before comparison with the tolerances.
Split conformal calibration chooses this factor so that both true errors are
no larger than the adjusted estimates with probability at least $1-\alpha$
under exchangeability (for example, independent cases from the same
distribution) \cite{Lei2018}; Appendix~\ref{app:technical-implementation}
gives the calculation. Thus, for a new case, an unsafe selection has probability
at most $\alpha$. This does not bound the unsafe fraction among accepted cases
or protect against distribution shift.

\papermethodbarrier
\section{Numerical results}\label{sec:experiments}

We evaluate the selector on independent cases from three benchmarks. Online
totals include feature construction, inference, the selected solve, and
fallbacks to the full model;
paired labeling and training are offline. Sampling, data splits, architectures,
discretizations, solver settings, refinement checks, and timing protocols are
collected in Appendix~\ref{app:technical-implementation}, after the references.
The code is available in the
\href{https://github.com/danielfdzm/learning-boundary-condition}{project repository}.

\subsection{Experiment 1: harmonic galvanic corrosion model with a boundary jump}
\label{sec:electrochem-setup}

Following \cite{Fernandez2026SingularLimit}, let
$\Omega=(0,0.02)\times(0,0.01)$. The bottom boundary is split into cathodic
and anodic halves, denoted by $\Gamma_c$ and $\Gamma_a$; the remainder
$\Gamma_N$ is insulated. The full problem is
\begin{equation}
 \begin{cases}
  -\Delta\phi_\kappa=0 & \text{in }\Omega,\\
  \partial_n\phi_\kappa=-\frac{1}{\kappa}i_c(\phi_\kappa)
    & \text{on }\Gamma_c,\\
  \partial_n\phi_\kappa=-\frac{1}{\kappa}i_a(\phi_\kappa)
    & \text{on }\Gamma_a,\\
  \partial_n\phi_\kappa=0 & \text{on }\Gamma_N,
 \end{cases}
 \label{eq:electrochem-full-model}
\end{equation}
where the cathodic and anodic currents are
\begin{align}
 i_c(\phi)
 &=i_{c,0}\left[
   \exp\!\bigl(C_1(\phi-\phi_c)\bigr)
   -\exp\!\bigl(-C_2(\phi-\phi_c)\bigr)\right],\nonumber\\
 i_a(\phi)
 &=i_{a,0}\left[
   \exp\!\bigl(A_2(\phi-\phi_a)\bigr)
   -\exp\!\bigl(-A_1(\phi-\phi_a)\bigr)\right].
 \label{eq:electrochem-currents}
\end{align}
We use the fixed exponential slopes from
\cite{Fernandez2026SingularLimit} and vary $\kappa$, the equilibrium
potentials, and the current scales around their reference values. The varied
parameter vector is
$(\kappa,\phi_a,\phi_c,i_{c,0},i_{a,0})$. The limit sets
$\phi_0=\phi_c$ on $\Gamma_c$ and $\phi_0=\phi_a$ on $\Gamma_a$, so its trace
jumps at their junction. Both models use the same graded conforming $P_1$
mesh. Errors are measured in $\Omega$ and on
$\Gamma_\star=\Gamma_c\cup\Gamma_a$. The estimator corrects two linearized
error indicators rather than learning the errors from scratch.
In each repetition, split conformal calibration uses 90 cases that are
separate from the training and evaluation sets.

At 5\% domain and boundary tolerances, the raw selector accepts most of the
cases that the paired reference identifies as safe. Split conformal calibration
reduces unsafe choices but misses more safe cases. Across
repetitions, the conditional unsafe rate ranges from 0.87\% to 3.48\% for the
raw selector and from 0 to 0.93\% after calibration.

Table~\ref{tab:butler-volmer-results} compares the neural rules with a tuned
$\kappa$ threshold, the two linearized indicators
used directly, and ridge residual regression using the same inputs and target.
The paired reference row is an unattainable oracle bound.

\paperbegincorrosiontable
\centering
\scriptsize
\setlength{\tabcolsep}{4pt}
\begin{tabular}{@{}lccc@{}}
\toprule
Estimator & Limit use [\%] & Unsafe / uses & Missed / safe \\
\midrule
Tuned $\kappa$ threshold & 33.9 & 20/542 & 48/570 \\
Linearized indicator & 0.5 & 0/8 & 562/570 \\
Ridge residual regression & 36.5 & 14/584 & 0/570 \\
Neural residual regression & 35.2 & 8/564 & 14/570 \\
Calibrated neural regression & 33.4 & 1/534 & 37/570 \\
Paired reference & 35.6 & 0/570 & 0/570 \\
\bottomrule
\end{tabular}
\tablecaptiongap
\caption{Galvanic corrosion results with our framework.}
\label{tab:butler-volmer-results}
\paperendcorrosiontable

Figure~\ref{fig:butler-volmer-decision} resolves these counts on the first
locked evaluation. The estimates follow the paired errors over nearly three
decades, including the stiff cases whose errors reach $80\%$, and the induced
decision differs from the paired reference on three of the 320 cases: one
unsafe selection and two safe cases sent to the full model. The conformal
factor of about $1.09$ withdraws four acceptances; this removes the unsafe
selection and raises the number of missed safe cases from two to five, which
is the trade-off reported above. Both error types are comparable in magnitude
for this problem, so neither tolerance is uniformly the binding one, and cases
near the corner of the acceptance region decide the outcome.

\paperbegincorrosionfigure
\centering
\includegraphics[width=\linewidth]{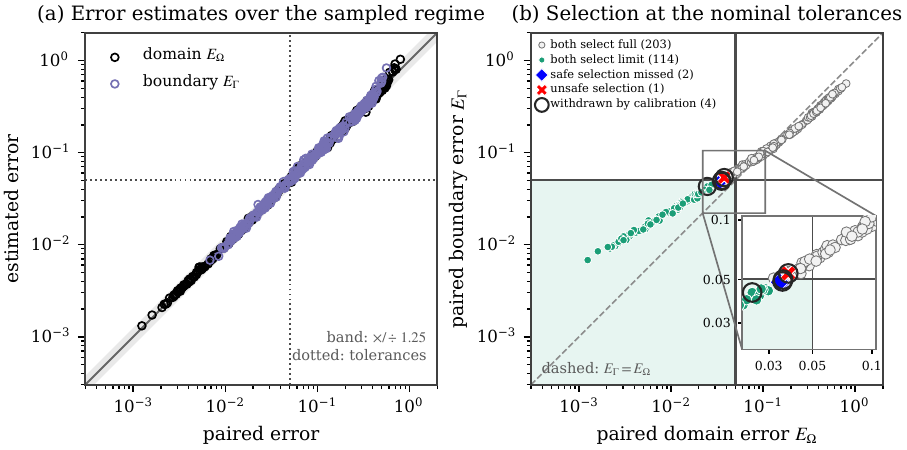}
\caption{Estimated errors and the resulting selection on the first locked
evaluation. (a) Estimated against paired errors, with the common
5\% tolerance dotted. (b) The same cases in the plane of paired errors; the
shaded rectangle is the acceptance region, the inset magnifies its corner,
and black rings mark the acceptances withdrawn by calibration.}
\label{fig:butler-volmer-decision}
\paperendcorrosionfigure

\papercorrosionbarrier

\subsection{Experiment 2: stationary nonlinear transfer}
\label{sec:stationary-high-speedup}

We consider the stationary problem on $\Omega=(0,1)^2$:
\begin{equation}
 \begin{aligned}
 -\Delta u_\kappa+u_\kappa&=f_\mu
   &&\text{in }\Omega,\\
 \partial_n u_\kappa+\frac{1}{\kappa}
 \bigl[(u_\kappa-g_\mu)+\gamma(u_\kappa-g_\mu)^3\bigr]&=0
   &&\text{on }\partial\Omega.
 \end{aligned}
 \label{eq:stationary-high-gap-full}
\end{equation}
Here $g_\mu=g_0+g_x\cos(2\pi x)+g_y\sin(\pi y)$. Let
$s_1=\sin(\pi x)\sin(\pi y)$ and $s_2=\sin(2\pi x)\sin(\pi y)$, so
$f_\mu=f_1s_1+f_2s_2$. We vary the stiffness, nonlinearity, boundary target,
and load amplitudes. The corresponding parameter vector is
$\mu=(\kappa,\gamma,g_0,g_x,g_y,f_1,f_2)$.

At the nominal tolerances, the selector misses one safe case. For this finite
test set, the descriptive 95\% upper endpoint for unsafe use is 6.2\%, so it is
not a guarantee.

\paperbeginstationarytable
\centering
{\small
\begin{tabular}{lrr}
\toprule
Quantity & Value & Unit \\
\midrule
Full solve, median & 130.597 & ms/query \\
NN + limit solve, median & 0.978 & ms/query \\
Accepted path, ratio of medians & 133.5 & $\times$ \\
Policy speedup & 8.3 & $\times$ \\
Raw selector limit choices & 58/64 & cases \\
Unsafe limit choices & 0/58 & cases \\
\bottomrule
\end{tabular}
}
\tablecaptiongap
\caption{Stationary selector accuracy and timings on one thread at nominal
tolerances. Policy speedup includes feature construction, inference, selected
solves, and all fallbacks to the full model.}
\label{tab:stationary-high-speedup}
\end{table}

Table~\ref{tab:stationary-high-speedup} shows a large acceleration on accepted
limit paths. Fallbacks reduce, but do not eliminate, the policy gain.

\FloatBarrier

\subsection{Experiment 3: evolutionary nonlinear transfer}
\label{sec:evolutionary-high-speedup}

The evolutionary benchmark uses the same square domain over $0<t\le0.8$:
\begin{equation}
 \begin{aligned}
 \partial_tu_\kappa-\nu\Delta u_\kappa&=f_\mu(x,y,t)
   &&\text{in }\Omega,\\
 \nu\partial_nu_\kappa+\frac{1}{\kappa}
 \frac{\sinh\!\bigl(\beta(u_\kappa-g_\mu)\bigr)}{\beta}&=0
   &&\text{on }\partial\Omega.
 \end{aligned}
 \label{eq:evolutionary-high-gap-full}
\end{equation}
The limit imposes $u_0=g_\mu$ on the boundary. We set
$g_\mu=0.65+a_g\sin(2\pi\omega_gt+\psi)+0.075(x-0.5)-0.055(y-0.5)$,
$f_\mu=a_f\cos(1.4\pi t)\exp[-35((x-0.68)^2+(y-0.34)^2)]$, and
$u(\cdot,0)=g_\mu(\cdot,0)+0.12\sin(\pi x)\sin(\pi y)$. We vary the stiffness,
diffusivity, nonlinearity, boundary oscillation, source amplitude, and phase.
The corresponding parameter vector is
$\mu=(\kappa,\nu,\beta,a_g,\omega_g,a_f,\psi)$.

At the nominal tolerances, the larger test set has one unsafe accepted case and
no missed safe case; the separate direct policy evaluation also has one
borderline boundary miss. Accepted limit paths are much faster than full solves,
but fallbacks substantially reduce policy speedup. The smaller policy evaluation
has correspondingly wide uncertainty, so both results are descriptive rather
than safety guarantees.

\FloatBarrier

\subsection{Tolerance versus speedup}
\label{sec:tolerance-speedup}

We hold each network fixed and multiply its tolerance pair by
$\lambda$. The tolerance pairs are $(0.5\%,0.5\%)$ for the stationary problem
and $(1\%,2\%)$ for the evolutionary problem. For the accepted set $A_\lambda$
defined by predictions, policy speedup is
\begin{equation}
 S_{\mathrm{policy}}(\lambda)
 =
 \frac{\sum_i T_{\mathrm{full},i}}
 {\sum_{i\in A_\lambda}T_{\mathrm{NN+limit},i}
  +\sum_{i\notin A_\lambda}T_{\mathrm{NN+full},i}}.
 \label{eq:tolerance-policy-speedup}
\end{equation}
Here $T_{\mathrm{full}}$ is the cost of the full solve, while
$T_{\mathrm{NN+limit}}$ and $T_{\mathrm{NN+full}}$ are the measured costs of
the selected branches. The denominator therefore includes prediction, limit
solves, and all full fallbacks.

\begin{figure}[htbp]
\centering
\includegraphics[width=0.88\linewidth]{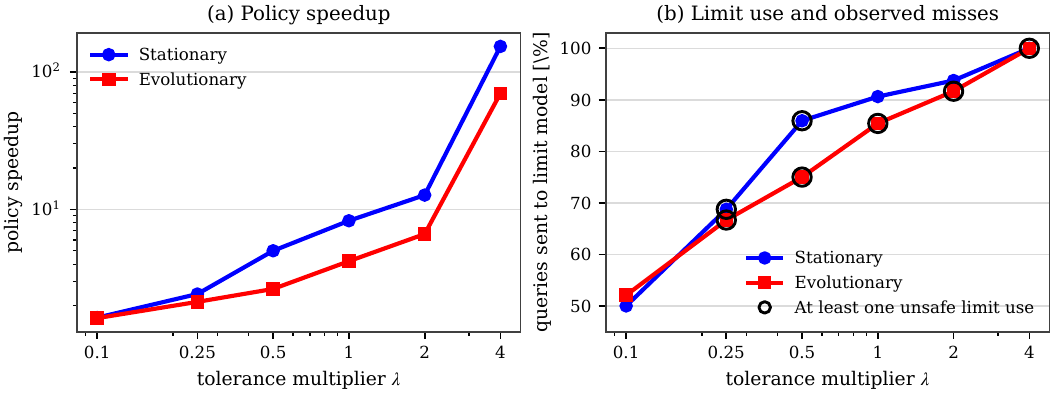}
\caption{Policy speedup and limit use versus tolerance multiplier
$\lambda$; open rings mark observed unsafe selections.}
\label{fig:tolerance-speedup-tradeoff}
\end{figure}

Relaxing the tolerances increases limit use and policy speedup in both
benchmarks. At $\lambda=1$, the stationary sweep agrees with
Table~\ref{tab:stationary-high-speedup}; the evolutionary result uses the larger
test set rather than the separate policy evaluation. Ringed markers denote an
observed unsafe selection, so speedup must be interpreted with the tolerances
and safety result.

\FloatBarrier

\section{Conclusion and open problems}\label{sec:conclusion}

We presented a framework for learning when one boundary condition can replace
another, focusing on singular limits from Robin or nonlinear laws to Dirichlet
data. Unlike full field neural surrogates, which are unreliable near singular
behavior, it predicts only domain and boundary errors. A PDE solver enforces the
choice, the full law provides a fallback, and tolerances can change without
retraining. Tests on stationary and evolutionary PDEs show that tolerance,
limit use, reusable operators, and fallbacks determine policy speedup.

This leaves several questions for the further development of the framework:
\begin{itemize}
  \item \textbf{Reliability among accepted cases.}
  The optional calibration bounds the marginal probability of an unsafe
  decision, but does not provide the same bound for a tolerance violation
  conditional on selecting the limit model. An open task is to control this conditional risk
  while retaining a useful acceptance rate, and to quantify the calibration
  sample size needed for that trade-off.

  \item \textbf{Accuracy throughout a trajectory.}
  The evolutionary test selects one model for an entire trajectory and measures
  errors at stored time levels. Controlling errors between these levels
  requires temporal error estimates. Adaptive switching during a solve also
  requires tracking accumulated error: returning to the full model does not
  by itself remove errors inherited from earlier limit solves.

  \item \textbf{Asymptotic features beyond linearization.}
  The corrosion estimator already corrects two linearized indicators.
  Extending this idea to other boundary laws could use justified convergence
  rates, boundary-layer structure, and junction corrections to improve
  predictions where data are sparse. Such features must account for changes
  in the rate or relevant norm caused by nonsmooth data.
\end{itemize}
\FloatBarrier

\section*{Acknowledgements}

The authors sincerely thank Enrique Zuazua for his unwavering support. This work was
funded by Schaeffler Technologies AG \& Co. KG.

\providecommand{\bysame}{\leavevmode\hbox to3em{\hrulefill}\thinspace}
\providecommand{\MR}{\relax\ifhmode\unskip\space\fi MR }
\providecommand{\MRhref}[2]{%
  \href{http://www.ams.org/mathscinet-getitem?mr=#1}{#2}
}
\providecommand{\href}[2]{#2}

\appendix
\numberwithin{table}{section}
\section{Technical implementation details}\label{app:technical-implementation}

\subsection*{\textbf{Common protocol}}
Held out evaluation and validation data do not split trajectories; fitting data
set all standardizations, and fixed seeds determine every design and fit
(Table~\ref{tab:experiment-protocols}).
For calibration, sort $s_i=\max_j(E_{i,j}/\widehat E_{i,j})$ on
$n_{\rm cal}$ cases unused for training. With
$k=\lceil(n_{\rm cal}+1)(1-\alpha)\rceil$, take the $k$th smallest score $q$
(or $+\infty$ if $k>n_{\rm cal}$) and set $c=\max(1,q)$.

\begin{table}[htbp]
\centering
\scriptsize
\setlength{\tabcolsep}{2.5pt}
\begin{tabular}{@{}L{0.12\linewidth}L{0.14\linewidth}L{0.18\linewidth}L{0.20\linewidth}L{0.26\linewidth}@{}}
\toprule
Study & Fit/cal./test & Online inputs & Discretization & Estimator fit \\
\midrule
Corrosion & $270/90/320$, five repeats & 8: five parameters, jump, two indicators & graded $31^2$ $P_1$ & $(96,96)$ ReLU/Adam; $10^{-4}$ penalty; $10^{-3}$ rate; 15\% early stop \\
Stationary & $256/\text{--}/64$ & 8 parameter/load features & $97^2$ $Q_1$ & $(24,12)$ tanh with L-BFGS; $10^{-4}$ penalty \\
Evolutionary & $128/32/(12+48)$ & 10 parameter features & labels: $25^2$/48 steps; checks: $97^2$/192 steps & $(32,16)$ tanh with L-BFGS; $2\times10^{-3}$ penalty \\
\bottomrule
\end{tabular}
\tablecaptiongap
\caption{Shared protocol; ``cal.'' means calibration or validation; test sets
are disjoint.}
\label{tab:experiment-protocols}
\end{table}

\subsection*{\textbf{Experiment 1: Corrosion}}
Ranges are $\log_{10}\kappa\in[-7,-3]$,
$\phi_a\in[-0.26,-0.14]$, $\phi_c\in[0.14,0.26]$,
$i_{c,0}\in[1.5,6]\times10^{-4}$, and
$i_{a,0}\in[1.5,6]\times10^{-2}$; currents are sampled logarithmically.
Calibration and test cases use independent samples from this design. Inputs are
the five parameters, the jump $\phi_c-\phi_a$, and two
linearized indicators. From the limit flux $q_h$, the correction
$\delta_\Gamma=-\kappa q_h/i_*'(\phi_*)$ and its harmonic extension give
$b_\Omega,b_\Gamma$; targets are $m_j-\log_{10}b_j$.

\subsection*{\textbf{Experiment 2: Stationary}}
Fitting and test samples use
$\log_{10}\kappa\in[-5,-0.5]$, $\log_{10}\gamma\in[4,8]$,
$g_0\in[0.8,1.2]$, $g_x\in[-0.3,0.3]$, $g_y\in[-0.25,0.25]$,
$f_1\in[0,8]$, and $f_2\in[-4,4]$. Inputs are $\log_{10}\kappa$,
$\log_{10}(1+\gamma)$, the other five values, and $\log_{10}(\kappa L)$, with
$L=1+|f_1|+|f_2|+4\pi^2(|g_x|+|g_y|)$.

\subsection*{\textbf{Experiment 3: Evolutionary}}
Samples use $\log_{10}\kappa\in[-3.5,-0.5]$,
$\nu\in[0.035,0.14]$, $\beta\in[20,120]$, $a_g\in[0.06,0.30]$,
$\omega_g\in[0.5,2.3]$, $a_f\in[0.15,1.10]$, and $\psi\in[0,2\pi]$.
Inputs are $\log_{10}\kappa,\allowbreak\nu,\allowbreak\beta/100,\allowbreak
a_g,\allowbreak\omega_g,\allowbreak a_f,\allowbreak\sin\psi,\allowbreak
\cos\psi,\allowbreak\log_{10}(\kappa/\nu)$, and $\beta a_g$. Labels use backward Euler and
checks use centered differences of order two. Full solves use damped Newton;
selected trajectories reuse the limit factorization.

\subsection*{\textbf{Timing and refinement}}
Single thread timings exclude training and the initial factorization.
Table~\ref{tab:stationary-high-speedup} uses 24 cases and seven repeats; policy
speedup uses 64 tests, three repeats, and six fallbacks. Sweeps time each path
once with prediction and fallback. Refining one case per study to $41^2$,
$129^2$, and $129^2$ (256 steps) changed errors by at most
$2.4\times10^{-3}$, $1.8\times10^{-7}$, and $2.2\times10^{-4}$; no decision
changed.

\FloatBarrier

\end{document}